\documentclass{article}

\begin{document}
\large
\title{Constants of Derivation and Differential Ideals}         

\author{Eloise Hamann\\Department of Mathematics,\\San Jose State University,San Jose, Ca 95192-0103\\hamann@math.sjsu.edu}
 
\maketitle
\textbf{Abstract:} Let R be a differential domain finitely generated over a differential field F of characteristic 0. Let C be the subfield of differential constants of F. This paper investigates conditions on differential ideals of R that are necessary or sufficient to guarantee that C is also  the set of constants of differentiation of the quotient field, E, of R.  In particular when C is algebraically closed and R has a finite number of height one differential prime ideals, there are no new constants in E. An example where F is infinitely generated over C shows the converse is false. If F is finitely generated over C and  R is a polynomial ring  over F, sufficient conditions on F are given so that no new constants in E does imply only finitely many height one prime differential ideals in R. In particular  F can be \={Q}(T) where T is a finite transcendence set.  \\

\begin{center}
\textbf{Introduction}
\end{center}

We fix notation for the paper.\\

\noindent\fbox{\begin{minipage}{12.5cm}
R = F[x$_{1}$,...,x$_{n}$] is a differential domain over the differential field F.\\ 
\noindent C = the characteristic 0 subfield of constants of differentiation of F.\\
\noindent E = the quotient field of R.\\
\noindent Q = the copy of the field of rationals $\subseteq$ C.\\
\noindent D denotes the derivation map.\\
\noindent\textbf{0.1}
\end{minipage}}\\

It is known  when C is algebraically closed that if R has no differential ideals then there are no new constants of differentiation in E. The converse is false. However, the question is raised  whether given no new constants in E, it is always possible to find a finitely generated differential F-subalgebra T of E containing R such that T has no proper differential ideals. [M] In [H1] there is an example to show that this need be the case.  The example is constructed with F of infinite transcendence degree over C.  It is natural to reformulate the question with the assumption that F be finitely generated over C.  [H1] purports to give a counterexample to this question but the example is incorrect. The author's work in calculating constants in various examples has led to the investigation of the relation between new constants and height one prime differential ideals.  Any  relation of new constants to other aspects of the ring is of interest because there is a theory of the nature of differential field extensions over a differential base field where the field extension was constructed to contain solutions to differential equations.  The theory requires no new constants for satisfying results. \\ 

 We assume basic familiarity with differential ideals (ideals closed under differentiation), differential homomorphisms (those that commute with differentiation) whose kernels are differential ideals,  the correspondence theorem for differential rings which are homomorphic images of a differential ring, differential ideals localize to differential ideals, and the fact that minimal primes over differential ideals are differential ideals in characteristic 0.  \\  

\section{A sufficient condition for no new constants}

We assume the basic setup of 0.1 of the introduction. \\  

\textbf{Lemma 1.1:}  If r $\in$ R is such that Dr = 0, then either r$-$ c is a nonunit in R for infinitely many c $\in$ C or r is algebraic over F.\\

\textbf{Pf:} Let Dr = 0.  The proof of this is essentially that of Lemma 1.16 of [M].  Let \={R} = R$\otimes$$_{F}$\={F} where \={F} is the algebraic closure of F.  If r $-$ c is a nonunit in \={R}, then clearly it's a nonunit in R.  View r as a function on the variety V defined by R.  r(V) in \={F} is constructible and therefore contains a dense open subset U of the Zariski closure of r(V).  If r(V) is infinite, then so is U so r(V) is cofinite in \={F}.  Thus, there are only a finite number of values c such that r(v) = c has no solution so r $-$ c vanishes for an infinite number of c and r $-$ c  is a nonunit for an infinite number of c.  If r(V) is finite, then it is a singleton, since R a domain implies V is irreducible and r is an \={F}-constant.  \\

\textbf{Theorem 1.2:}  Let R = F[x$_{1}$,...,x$_{n}$] be a differential domain over the differential field F of characteristic 0.  Let C be the algebraically closed field of constants of F.  If R has only finitely many height one differential prime ideals, then the quotient field E of R has no new constants of differentiation.\\ 

\textbf{Pf:} We prove the contrapositive, so suppose there exist f/g $\in$ E $-$ C such that D(f/g) = 0.  Apply Lemma 1.1 to f/g in R[1/g].  First suppose that f/g is algebraic over F and p(X)= $\sum_{i=0}^m$t$_{i}$X$^{i}$ is f/g's irreducible polynomial over F.\\

\noindent  Then $\sum_{i=0}^m$t$_{i}$(f/g)$^{i}$ = 0 implies that $\sum_{i=0}^m$Dt$_{i}$(f/g)$^{i}$ = 0. \\

  Thus, $\sum_{i=0}^m$Dt$_{i}$X$^{i}$ must be an F-constant multiple of p(X) so  Dt$_{i}$ = st$_{i}$ for some s $\in$ F.  Fix j so t$_{j}$ $\ne$ 0.  Then, for any other nonzero t$_{i}$ we have t$_{i}$Dt$_{j}$ = t$_{j}$Dt$_{i}$ so t$_{i}$/t$_{j}$ is a differential constant.  t$_{i}$/t$_{j}$ $\in$ F so t$_{i}$/t$_{j}$ $\in$ C.  Let t$_{i}$ = c$_{i}$t$_{j}$.  Then $\sum_{i=0}^m$c$_{i}$(f/g)$^{i}$ = 0.  However, C is algebraically closed so p(X) must be linear and therefore f/g $\in$ C contradicting our choice of f/g.  Thus, f/g is not algebraic over F and by 1.1 we must have that f/g $-$ c is a nonunit in R[1/g] for an infinite number of c $\in$ C.  Since (f/g $-$ c) is a proper differential ideal  infinitely often we obtain an infinite number of height one primes P$_{c}$ minimal over (f/g $-$ c).  Then Q$_{c}$ = P$_{c}$ $\cap$ R is a differential height one prime ideal of R.  Since R[1/g] is a localization of R, the Q$_{c}$ are an infinite collection of distinct height one differential primes of R.\\ 

While it is customary to assume that C is algebraically closed in working with differential algebras, milder assumptions suffice for some purposes.  The example R = Q[x,y] $\cong$ Q[X,Y]/P with P = (X$^{2}$ + Y$^{2}$), DX = X, DY = Y, and Q is the rationals can be shown to have (x,y) as its only height one prime differential ideal but x/y is a differential constant of R not in Q. (Work in Q[X,Y] and pass to \={Q}[X,Y] to get (X,Y) the only maximal ideal of Q[X,Y]containing P and differential.) Thus, the assumption that C be algebraically closed is necessary.

\section{Necessary condtions for no new constants}

We continue to assume the setup in 0.1 of the introduction.  In addition we assume that F is finitely generated as a field over C.  The reason for this assumption is illustrated by the following example.\\  
 
\textbf{Example:}  Let F = C(\{$\alpha_{i}$\}), i $\in$ Z$^{+}$ where $\alpha_{i}$ are independent indeterminates over any field C of characteristic 0.  Define a derivation D on F by Dc = 0 if c $\in$ C and D($\alpha_{i}$) = $\alpha_{i}^{3} - 2\alpha_{i}^{2} + 2\alpha_{i}$ and extend to F by the rules of derivation.  Let R = F[X], E the quotient field of R and extend D to E by defining D(X) = X$^{3} - 2X^{2} +2X$.\\

In [H1], it is shown that C = the constants of differentiation of E.  It is easy to see that each (X $-$ $\alpha_{i}$) is a differential height one prime ideal of R and there are an infinite number.  Thus, the condition that R have a finite number of height one prime differential primes is not necessary in order for there to be no new constants in E. The situation in this example, however, is related to the fact that F is not finitely generated over C.\\

The following lemma is primarily a useful observation.\\

\textbf{Lemma 2.1}  Let R be a differential domain over a differential field F, with derivation D.  Let D$'$ = fD where f $\ne$ 0 $\in$ F.  Then an ideal I of R is a differential ideal under D iff I is a differential ideal under D$'$. Further w is a D-constant iff w is a D$'$ constant.\\  

\textbf{Pf:} Since 1/f $\in$ F, it suffices to prove one direction.  Let I be a differential ideal of R under D.  Then let r $\in$ I.  D$'$r = fDr $\in$ I since Dr $\in$ I. Since we are working in a domain the claim about constants is clear.\\   

\textbf{Lemma 2.2}  If R = F[x$_{1}$,...,x$_{n}$] is any differential domain over a field F and R is integral over a differential subdomain A with only finitely many height one prime differential ideals then the same holds for R. In particular, if F is algebraic over a subfield K such that A = K[x$_{1}$,...,x$_{n}$] is a differential domain with only finitely many height one prime differential ideals, then R has only finitely many.\\

\textbf{Pf:} If P is a height one prime differential ideal of R, Q = P$\cap$A is such an ideal of A.  Since R is Noetherian and there are only finitely many possibilities for Q, R can have only finitely many height one prime differential ideals, namely those minimal over extensions of those of A. If F is algebraic over K, then R will be integral over K[x$_{1}$,...,x$_{n}$]. \\

We also freely use the fact that if a differential domain has only finitely many height one prime differential ideals so does any localization.  Localization by a multiplicative set S is easily seen to preserve the differential property of an ideal and there is a one to one order preserving correspondence between the primes of R disjoint from S and those in R$_{S}$ so the number of differential primes of R$_{S}$ of a given height is $\leq$ the number in R.\\
  
The following list of results are nearly cumulative, but are listed separately for clarity.\\

\textbf{Proposition 2.3:}  Let R be the differential polynomial ring C[X$_{1}$,...,X$_{n}$] where C is any subfield of  \={Q},  the algebraic closure of Q.  If E = q.f.(R) has no new constants, then R has only finitely many height one prime differential ideals.\\
 
\textbf{Pf:}  Let S be the finite set of coefficients required to express each DX$_{i}$ as polynomials in R.  Thus, A = Q[S][X$_{1}$,...,X$_{n}$] is a differential subring of R. Clearly, there are no new constants in q.f.(A) so by Lemma 2 of [H1], A has only finitely many height one differential prime ideals. The result follows from 2.2.\\

\textbf{Proposition 2.4:} Let R be the differential polynomial ring  C[X$_{1}$,...,X$_{n}$] where C is the set of differential constants of R and C = C$'$(T) where C$'$ is a subfield of \={Q} and T is a pure transcendence set of any cardinality. If E = q.f.(R) has no new constants, then R has only finitely many height one prime differential ideals.\\

\textbf{Pf:} Let S be a finite subset of C$'$ and  T* be a finite subset of T such that \{DX$_{i}$\} $\subseteq$ A = Q[S](T*)[X$_{1}$,...,X$_{n}$] so that A is a differential subring of R.  By 2.1 we may replace D with a constant multiple  so that B = Q[S][T*][X$_{1}$,...,X$_{n}$] is a differential  subdomain with no change in the differential ideals of R. We show that B has only finitely many height one prime differential ideals (y) with Dy = zy and z $\ne$ 0.\\
Let M be the maximum total degree that occurs in any DX$_{i}$ in the variables t$_{i}$ $\in$ T and X$_{i}$. Now let V be the vector space over Q generated by products of elements of a vector space basis for Q[S] with monomials in t$_{i}$ and X$_{i}$  of degree at most M.    V is a finite dimensional vector space.  Let Dw = zw.  Dw = ($\sum\partial$w/$\partial$ X$_{i}$) DX$_{i}$.  Each partial derivative involves terms of no higher degree than that of w in X$_{i}$ and t$_{i}$.  Then the degree of ($\partial$w/$\partial$ X$_{i}$)DX$_{i}$ cannot exceed the degree of w by more than  M in the variables X$_{i}$ and t$_{i}$.  Thus,  z $\in$ V.\\
Let L = \{v $\in$ V such that there exists w in B with Dw = vw\}.  Since V is a finite dimensional vector space over Q we can choose a maximal finite linearly independent subset L* =  \{z$_{i}$\} of L  associated with prime elements w$_{i}$ such that Dw$_{i}$ = z$_{i}$w$_{i}$ with z$_{i}$ $\ne$ 0.  Let (y) be a nonzero prime differential ideal with  Dy = zy $\ne$ 0, we must have z in L so z and L* are linearly dependent over Q.  Let nz = $\sum$n$_{i}$z$_{i}$ where n $\in$ Z$^{+}$, n$_{i}$ $\in$ Z.  Let b = $\prod$w$_{i}^{n_i}$ so Db = $\sum$n$_{i}$z$_{i}$b = nzb.  Also, Dy$^{n}$  = nzy$^{n}$ so D(y$^{n}$/b) = 0 and y$^{n}$/b is a differential constant.  Thus, y$^{n}$ = cb  for some c in C $\cap$ qf(B) = Q[S](T*). Let c = p/q with p,q $\in$ Q[S][T*] so qy$^{n}\prod$w$_{j}$$^{-n_{j}}$ = p$\prod$w$_{k}^{n_{k}}$ where n$_{j}$ is negative, n$_{k}$ is positive. Since (y) is a prime ideal, p $\in$ (y) or some w$_{k}$ $\in$ (y).  Since Dy $\ne$ 0, y $\notin$ Q[S][T*] which implies y has positive degree in the X$_{i}$ so some (w$_{k}$) must = (y).    Thus, the set of principal prime differential ideals not generated by constants is the finite set \{(w$_{i}$)\}.\\  
Since A is a localization of B, A's height one prime differential ideals are extended from those of B.  Since the constants of differentiation in A are invertible there are no proper ideals (y)  with Dy = 0.  Thus, A's height one prime differential ideals are the finitely many (w$_{i}$)A.  Let T** = T $-$ T*.   Clearly T** is an algebraically independent set over A and A[T**] is a differential ring under the current derivation since we are only adding constants of derivation.  We argue that A[T**] has only finitely many height one prime differential ideals.    Let (y) be such an ideal and let y = $\sum_{i=0}^m$a$_{i}$m$_{i}$ where a$_{i}$ $\in$ A and m$_{i}$ is a monomial in elements of T**.  Then Dy = $\sum_{i=0}^m$Da$_{i}$m$_{i}$.  Since Da$_{i}$ $\in$ A, we must have Da$_{i}$ = za$_{i}$ with z $\in$ q.f.(A).  Thus, Da$_{i}$/a$_{i}$ = Da$_{j}$/a$_{j}$, a$_{i}$Da$_{j}$ $-$ a$_{j}$Da$_{i}$ = 0 and a$_{i}$/a$_{j}$ is a differential constant. Since E has no new constants, a$_{i}$ = c$_{i}$a$_{j}$ for a$_{j}$ a fixed nonzero coefficient and c$_{i}$ $\in$ C$\cap$q.f.(A) = C$'$(T*).  Then y = ($\sum_{i=0}^m$c$_{i}$m$_{i}$)a$_{j}$ and  (y) = (a$_{j}$) in A[T**].  Clearly  (a$_{j}$) in A is a height one prime differential ideal so (y) is merely an extension of one of the finite number belonging to A.  Localizing we obtain Q[S](T)[X$_{1}$,...,X$_{n}$] has only finitely many height one prime differential ideals.  Finally, R is integral over Q[S](T)[X$_{1}$,...,X$_{n}$] since we obtain R  by adjoining elements algebraic over Q[S] so R also has only finitely many height one differential prime ideals by 2.2.\\

\textbf{Theorem 2.5:} Let R = F[X$_{1}$,...,X$_{n}$] be a differential polynomial ring where F is a differential field of characteristic 0 of the form C(W) where W is a finite transcendence set of F over C, the subfield of differential constants.  Assume C = C$'$(T) where T is a transcendence set of any cardinality over C$'$, any subfield of \={Q}.   If E = q.f.(R) has no new constants, R has only finitely many height one  prime differential ideals.\\  

\textbf{Pf:} There exists f $\in$ F so fDw$_{i}$ and fDX$_{i}$ are in A = C[W][X$_{1}$,...,X$_{n}$]for all i.   A is a differential subdomain under D$'$ = fD with no change in the differential ideals of R.  Since A is a finite polynomial ring over C with no new constants, 2.4 implies that A has only finitely many height one differential ideals.  Since R is a localization of A, the claim follows.\\

\textbf{Theorem 2.6:}  Let R = F[x$_{1}$,...,x$_{n}$] be a differential domain of Krull dimension d with  F = C(W), W a finite transcendence set over C, C = C$'$(T), T a transcendence set of any cardinality, C$'$ any subfield of \={Q}.    Suppose there are no new constants in E = q.f.(R). If R has a differential subring which is a polynomial ring in d variables over F, then R has only finitely many height one differential primes. In particular if R has a Noether normalization which is a differential subring, R has only finitely many height one prime differential ideals.   \\

\textbf{Pf:}   If R has a polynomial differential subring S in d variables over F, then R is algebraic over S.  Since R is finitely generated over S, there is a single element s in S such that R[1/s] is integral over S[1/s].  Only finitely many height one primes of R or S can contain s so R[1/s] (S[1/s]) has finitely many height one differential primes iff R (S) does. Since S is a polynomial ring over F, S has only finitely many height one prime differential ideals by 2.5.  Now apply 2.2 to the differential domain R[1/r].  \\

2.6 is only one example of how 2.2 might be used to obtain more general results.\\

The following result strengthens Theorem 3 of [H1].  \\  

\textbf{Theorem 2.7:}  Let R = F[X$_{1}$,X$_{2}$] be a differential polynomial ring over F  where   F = C(W), W a finite transcendence set over C, C = C$'$(T), T a transcendence set of any cardinality, C$'$ any subfield of \={Q}.   Suppose there are no new constants in E = q.f.(R).  Then R has a finitely generated differential extension within E with no nonzero proper differential ideals.\\

\textbf{Pf:} If a differential ring  has no nonzero prime differential ideals, then it has no nonzero proper differential ideals since it is known that a maximal differential proper ideal is prime.  Thus, suppose R does have nonzero prime differential ideals.   By 2.5, there are only finitely many height one differential prime ideals. DX$_{i}$ $\ne$ 0 since there are no new constants.  The height two ideals are maximal.  If M is maximal ideal which is differential there are finitely many maximal ideals of \={F}[X$_{1}$,X$_{2}$] lying over M\={F}[X$_{1}$,X$_{2}$] which  have  form (X$_{1}$ - c$_{1}$, X$_{2}$ - c$_{2}$) and are differential so all contain DX$_{i}$. Thus, DX$_{i}$ $\in$ M. If there are height one differential primes, let t be a nonzero element of the finite intersection.  Otherwise, let t = 1.  Then tDX$_{i}$ belongs to every prime differential ideal.  T = R[1/tDX$_{1}$] is a finitely generated differential extension of R.  Any nonzero prime differential ideal of T would intersect R in a nonzero prime differential ideal, but no nonzero prime differential ideal extends to a proper ideal in T so T has no nonzero prime differential ideals whence no nonzero proper differential ideals.     \\  

Algebraic geometers may be most interested in the case F = \={Q}(W).\\

\textbf{Concluding remarks:} While the connection between height one prime differential ideals and new constants could be further pursued, one can also ask about the connection between new constants and principal prime differential ideals.   No new constants requires  DY = zY for a fixed z to have only one solution for Y up to a constant multiple. When R is a UFD, this condition is sufficient. The condition relates well to principal differential ideals and perhaps the positive results in this paper about polynomial rings and height one prime ideals should be viewed in this vein.\\

\begin{center}
\textbf{Bibliography}
\end{center}

[H1]  E. Hamann, A Question on Differential Ideals,  \textit{Communications in Algebra}, Vol. 31, number 12 2003 6079-6091.

[H2]  E. Hamann, Corrigendum to A Question on Differential Ideals, to appear in \textit{Communications in Algebra}.
 
[K1] I. Kaplansky,An Introduction to Differential Algebra, Paris. Hermann 1957.

[K2] I. Kaplansky, Commutative Rings, Allyn and Bacon, Boston. 1970
  
[L] T. Little, Algebraic Theory of Differential Equations, Thesis San Jose State University, 1999.

[M] A. Magid, Lectures on Differential Galois Theory, University Lecture Series Vol. 7, Providence, AMS, 1991\\

\end{document}